# MATHEMATICS EDUCATION RESEARCH MUST BE USEFUL FOR THE CLASSROOM

Anthony A. Essien[1], Tanya Evans[2], Maitree Inprasitha[3], Salomé Martínez Salazar[4], Demetra Pitta-Pantazi[5]

University of the Witwatersrand [1], University of Auckland[2], Khon Kaen University[3], Universidad de Chile[4], University of Cyprus[5]

*The theme of the PME-48 conference, Making sure that mathematics education research reaches the classroom, implicitly acknowledges a concern: not all mathematics education research finds its way into classroom practice. This raises several fundamental questions: Which research fails to reach the classroom? Why? And should all research be expected to do so? It is therefore fitting that the plenary panel engages with these significant issues by debating the following motion: Mathematics education research must be useful for the classroom. This paper presents the debate as structured for the purposes of this publication. Following an introduction by the Chair, Tanya Evans, each panellist's arguments are presented individually. Anthony Essien and Salomé Martínez Salazar argue against the motion, while Maitree Inprasitha and Demetra Pitta-Pantazi argue in support. We hope this debate will stimulate ongoing dialogue and encourage the mathematics education research community to critically engage with this important issue—one that is central to the relevance and influence of research in our field.*

## INTRODUCTION (TANYA EVANS)

To debate the motion *"Mathematics education research must be useful for the classroom"*, it is first necessary to ground our discussion in the fundamental question: *What do we mean by mathematics education research?* Without a shared understanding of what mathematics education research entails, discussions about its usefulness—particularly to classroom practitioners—risk becoming either vague or polarising.

A useful frame for exploring this question is offered by Inglis and Foster (2018), who conducted a seminal study providing a comprehensive and systematic overview of mathematics education research over the past five decades. They analysed the full text of all articles published in two leading journals of the field—*Educational Studies in Mathematics* and the *Journal for Research in Mathematics Education*—from their inception to the present. Using a statistical modelling technique known as topic modelling, they identified prominent topics and traced their evolution over time. Importantly, Inglis and Foster drew on Lakatos's (1978) methodology of scientific research programmes to interpret their findings. The central idea in Lakatos's account is that the base descriptive unit of research is not, as Popper (1959) maintained, the individual hypothesis or even theory, but rather the research programme—a broader and more durable structure that guides the development of knowledge over time.

Lakatos (1978) described a research programme as comprising three main components: a "hard core" of foundational assumptions, a "protective belt" of auxiliary hypotheses, and a "heuristic" that governs how research progresses. The hard core defines the programme's core commitments and is shielded from falsification. The protective belt can be modified to accommodate anomalies from new findings, allowing the hard core to remain intact. The heuristic consists of guiding principles for selecting methods, investigating phenomena, generating data, and developing theories within the programme.

Importantly, Lakatos used this framework to distinguish between *progressing* and *degenerating* research programmes. A progressing programme generates new insights and explains more phenomena over time. A degenerating one, in contrast, responds to anomalies only through post hoc adjustments to its protective belt, without generating new predictive power. Over time, researchers may abandon a degenerating programme in favour of one that is progressing.

Inglis and Foster (2018) applied this framework to the history of mathematics education research, illustrating how various research programmes have risen and declined over time. For example, the constructivist research programme, which dominated the 1980s and 1990s, was founded on the hard-core assumption that knowledge construction is an individual process aimed at helping learners make sense of the world. According to this view, individuals receive sensory input, filter it, and actively organise it into mental schemas. Learning occurs when people construct new knowledge by creating new schemas or reorganisig existing ones. Rather than passively absorbing information, knowledge is actively constructed by the individual. The programme's heuristic was based on Piagetian clinical interviews and small-scale qualitative studies focused on how students make sense of mathematical ideas (Piaget & Cook, 1952).

By the late 1990s and early 2000s, a shift was underway: the so-called "social turn" in mathematics education research, documented by Lerman (2000), signalled the rise of a new sociocultural research programme. This programme rejected the individualist assumptions of constructivism and instead asserted that thinking and learning are fundamentally social processes. As Vygotsky (1986) famously put it, "the true direction of the development of thinking is not from the individual to the social, but from the social to the individual" (p. 36). Consequently, the heuristic shifted from individual interviews to analyses of classroom discourse, peer interaction, and the broader social and cultural contexts of mathematics learning (Inglis & Foster, 2018).

Rather than replacing the earlier paradigm entirely, the sociocultural turn ushered in a period of theoretical diversification. As Inglis and Foster (2018) observed, multiple research programmes—such as those based in semiotics, embodied cognition, and didactical theories—now coexist and often compete for attention. From a Lakatosian perspective, this theoretical pluralism can be seen as a vital indicator of disciplinary maturation. Different research programmes challenge each other by offering rival explanations and methodological innovations.

Crucially, the interaction between competing research programmes should enable researchers to uncover anomalies, refine protective belts, and drive the ongoing evolution and vitality of the field. However, as Dreyfus (2006) warned, theoretical proliferation without empirical challenge can also lead to fragmentation. He criticised the field for "inventing theories at a pace faster than we produce data to possibly refute our theories" (p. 78). In other words, the discipline risks hosting degenerating programmes that persist unchallenged—not because they are progressing, but because they are insulated from critique. To the extent that different research programmes fail to engage with one another, mathematics education research may lack the kind of productive tension that Lakatos considered essential for progress.

An illustrative example of the need for cross-programme dialogue lies in the relationship between the sociocultural research programme and the experimental psychology programme—so termed by Inglis and Foster—which is closely related to what is now often referred to as the *Science of Learning* (National Academies of Sciences, 2018)). The latter is characterised by a hard-core commitment to developing and testing theories through the manipulation of experimental stimuli to reveal underlying psychological processes—whether cognitive, behavioural, social, or developmental. Instructional approaches associated with this programme—such as explicit instruction, mastery learning, spaced retrieval practice, memorisation of arithmetic facts, and time-limited practice—have produced robust evidence of effectiveness (e.g., Fuchs et al., 2021; Gersten et al., 2009; Kirschner et al., 2006; McNeil et al., 2025). From a sociocultural perspective, however, these methods may appear to conflict with its hard-core assumptions, particularly the view that learning is constructed through social interaction and participation. Yet such findings constitute potential anomalies—ones that should prompt reconsideration of aspects of the protective belt of sociocultural theory (Inglis & Foster, 2018).

Conversely, experimental psychology might benefit from insights generated within sociocultural programme. For instance, research on micro-identities (e.g., Wood, 2013) shows that students' mathematical identities can shift rapidly in response to subtle contextual cues. This finding challenges a key assumption of many cognitive experiments: that individual traits and behaviours are relatively stable over the course of an experimental task. If micro-level identity shifts influence learning outcomes, then experimental psychology may need to refine its methods and assumptions to accommodate such variability (Inglis & Foster, 2018).

This brief account demonstrates that mathematics education research is not a monolithic endeavour but a diverse and evolving field comprising multiple research programmes with competing claims and heuristics. Against this backdrop, the question of usefulness—central to this plenary panel—becomes more complex. What counts as "useful" depends heavily on the epistemological commitments of the research programme in which the work is situated. For instance, a programme grounded in experimental psychology may value generalisable findings about cognitive processes, often derived from controlled studies with individual learners. By contrast, sociocultural research may prioritise rich, contextualised understandings of classroom interaction, student identity, or institutional discourse— insights that might not be readily transferable to other settings, but that can deepen our understanding of mathematics education more broadly. To ask both programmes to produce the same kind of "useful" knowledge is to ignore their differing assumptions, aims, and methods.

The theme of the PME-48 conference—Making sure that mathematics education research reaches the classroom—implicitly acknowledges a concern: that some research in mathematics education has limited impact on practice. But this prompts deeper questions. Which types of research struggle to reach classrooms, and why? Should all research be expected to do so? And perhaps most fundamentally, is usefulness to classroom practice a necessary criterion for mathematics education research to be considered valuable?

This plenary panel paper addresses these questions through a structured debate of the motion: *Mathematics education research must be useful for the classroom*. Following this introduction, the arguments of the four panellists are presented. Anthony A Essien and Salomé Martínez Salazar argue against the motion, while Maitree Inprasitha and Demetra Pitta-Pantazi argue in support. The perspectives they offer are shaped by the research traditions, methodological tools, and theoretical assumptions they have engaged with—elements that align in various ways with different research programmes as describe by Lakatos (1978). Framing the debate in these terms highlights how underlying assumptions about knowledge, methods, and aims influence what is considered "useful" in mathematics education research. Rather than seeking resolution, this discussion aims to illuminate the productive tensions across programmes and to stimulate further reflection on how diverse forms of research can contribute meaningfully to our field.

## PITFALLS OF CONSEQUENTIALISM IN MATHEMATICS EDUCATION (ANTHONY A ESSIEN)

### Mathematics and mathematics education research

Should mathematics be taught solely for its practical value to humanity? Should research in mathematics education be deemed useful only in as much as it informs classroom practice? An answer in the affirmative would suggests the consequentialist perspective of mathematics education and mathematics education research where the outcome or results are prized above the inherent benefit that accrue to a particular activity or act, in the present case, research in mathematics education. Speaking on the irrationality of mathematics education research, Herbert Wilf (2011) made a radical statement that "Research in mathematics education is generally hopelessly unintelligent. The field as a whole is not far from worthless." My own position is very different from Wilf's (2011) fundamentalist position. I argue rather against solely viewing mathematics education research though the consequentialist lens.

## The risks of consequentialism in mathematics education research

As its name suggests, consequentialism is a philosophical approach that prioritises the consequences or outcomes of actions when determining their value. It asserts that the overall consequences of an action encompass both the action itself and all that it brings about. Central to consequentialism is the idea that the outcomes of our actions possess inherent value and should serve as the main standard for value judgement. This view emphasises the importance of results above the methods used to achieve them.

One significant challenge with consequentialism, however, is that it does not specify which types of consequences are considered good. A parallel can be drawn with the statement: *Mathematics education research must be useful for the classroom*. What type of research does this imply, and by who? As an example, for a very long time before the early 1970s, research in mathematics education led the mathematics education community into believing that due to dissipation of stock of available intellect resulting from knowing more than one language, bi/multilingualism has a negative cognitive effect. For example, research by MacNamara in the early 1960s on the mathematics achievement of learners from English speaking homes when taught arithmetic in the medium of Irish falls into this category. MacNamara's sample consisted of 1,084 learners from 119 schools in Ireland. Each of the learners was given tests in problem arithmetic and mechanical arithmetic (see for example, Macnamara, 1977). The conclusion MacNamara drew from this study was that Irish bilingual education, and bilingualism in general, has a negative consequence (Baker, 1988). Today, with a plethora of further research in the field of language and mathematics (see Essien and Adler, 2025), the epistemic potential of the multiplicity of languages is acknowledged. Clearly, the work of MacNamara which used research results to argue for the compartmentalisation of human brain was deemed 'useful' in classroom practices, and in deciding on bi/multilingual policies at a time, albeit, with devastating consequences.

## Mathematics has a value in itself, so does mathematics education research

Mathematics possesses inherent worth, offering a distinct lens through which to understand the world by emphasising logic, abstract thinking, and the pursuit of knowledge independent of immediate practical use. In the words of Ernest (2016), "[I]rrespective of its worldly origins and applications, like many crafts and practices, mathematics takes on a life of its own". Mathematics engages the mind in a rich intellectual experience, nurturing critical thinking and problem-solving while encouraging the exploration of complex ideas. The process of engaging with mathematics can be deeply fulfilling, providing a sense of achievement and insight into the structures underlying various aspects of life. Furthermore, developing mathematical skills sharpens logical reasoning and analytical thinking—abilities that are valuable across academic disciplines, daily life, and professional fields.

By extension, mathematics education research is intrinsically valuable as it enhances our comprehension of how mathematical thinking develops and is best taught. While it serves practical purposes in refining teaching practices and curricula, its broader significance lies in promoting intellectual growth, fostering equitable learning environments, and reinforcing the role of mathematics as a vital component of human inquiry and progress.

**Stifling of intellectual exploration**

A consequentialist view of mathematics or mathematics education research can potentially stifle intellectual exploration, innovation and transdisciplinarity. While the value of research for classroom practice cannot be denied, placing a requirement that mathematics education research MUST be useful for the classroom undermines the broader purpose of scholarship and is a disregard for the multifarious dimensions inherent in mathematics education as a field. A consequence of this is that theoretical advancements are also stifled.

To the above also lies the question of immediate practicality vs long term applicability. The danger of consequentialism in mathematics education research lies in the fact that research would focus on the immediate mathematics education problems, and research on these. History has taught us that many breakthroughs in science and mathematics began with abstract investigations. Focusing on the immediate problem has the potential of stifling future innovations and in mathematics education research, this could translate into stifling critical research that could challenge dominant paradigms if these are not immediately applicable to classrooms. For instance, research on feminism in mathematics education may not offer immediate practical utility for the classroom, but could expose systemic inequalities.

**Conclusion**

In summary, mathematics is valuable not only for its practical applications but also for its intrinsic intellectual qualities, its role as a foundation for other disciplines, and its contribution to the development of essential cognitive skills and ethical values. In the same vein, while classrooms applicability is an important facet of mathematics education research, it cannot be the sole criterion of its value.

## SHOULD MATHEMATICS EDUCATION RESEARCH ALWAYS BE USEFUL FOR THE CLASSROOM? (SALOMÉ MARTÍNEZ SALAZAR)

In this paper, I argue that restricting research in mathematics education to approaches or areas deemed directly useful for classroom practice narrows the scope of inquiry and undermines disruptive academic efforts essential for advancing knowledge in any discipline. My arguments are, of course, shaped by my own academic trajectory. A significant part of my career has been dedicated to research in mathematics. From that perspective, it is difficult for me to imagine imposing on a discipline such a heavy burden as that of being 'useful.' Nevertheless, the discussion is far more nuanced, particularly given that mathematics education as a discipline is focused on the study of the teaching and learning of mathematics.

Indeed, several scholars have argued that education, as a discipline, should be directly oriented toward improving practice rather than primarily generating theoretical knowledge. This claim has been justified on the basis *of the distinctive character of human social action or, alternatively, by appealing to the kind of 'interest' that ought to guide educational—and perhaps all social—research* (Hammersley, 2003). Given the well-documented concern that mathematics educational research has had too little impact on improving classroom teaching and learning (Hiebert et al., 2002), it may seem tempting to restrict the focus of mathematics education research to that which provides a useful knowledge base for teaching. However, doing so entails significant risks.

### The necessity of freedom for research endeavors

Abraham Flexner, founder of the Institute for Advanced Study in Princeton, argued in his influential essay "The usefulness of useless knowledge" that many of the most significant scientific breakthroughs have originated from research unconcerned with immediate practical application. His central claim is that the freedom to pursue knowledge for its own sake is not a luxury, but a prerequisite for genuine intellectual progress. Although Flexner was likely referring to fields such as mathematics, physics, and the humanities, the spirit of his argument should extend to educational research as well. If we accept mathematics education as a scientific discipline grounded in the pursuit of knowledge, it too requires the freedom to explore foundational and theoretical questions aiming to understand the nature of mathematical thinking, teaching, and learning —regardless of their immediate applicability—if it is to contribute meaningfully to long-term understanding (Schoenfeld, 2000).

The unpredictable nature of knowledge development makes it unreasonable to constrain research questions, stimuli, or approaches to what is presumed to be 'useful' in the short term (Hammersley, 2003). Knowledge production is not a linear process; it requires open-ended exploration and can follow uncertain trajectories over extended periods (Flexner & Dijkgraaf, 2017).

**The need for interdisciplinary research**

Requiring short-term utility from research can restrict interdisciplinary approaches. According to Repko and Szostak (2020), for two research areas to be effectively integrated and produce new knowledge, genuine interdisciplinarity is essential: one that goes beyond merely combining perspectives, and instead involves the development of shared theoretical frameworks and clearly defined goals. This process demands ongoing dialogue that helps overcome linguistic, conceptual, and methodological barriers in order to establish a common language (Klein, 2017). It also requires epistemological and methodological flexibility, enabling researchers to adapt their approaches and open new avenues for exploration and integrated knowledge (Newell, 2001).

An example of the potential of interdisciplinarity for advancing toward a more comprehensive understanding of phenomena involved in mathematical learning is proposed by Leikin (2018). In her work, she addresses the potential advantages of connecting cognitive neuroscience with mathematics education, even though she acknowledges that they are still *tangent areas*, pointing out that *even though a relatively large number of neurocognitive studies have addressed numerical cognition, these studies are rooted in cognitive psychology and are not connected to the findings of mathematics education research.* Particularly, the use of technical terminology and the fact that many studies are performed in laboratory settings, which are not similar to classroom settings, make cognitive neuroscience research studies have little impact on the processes of learning primary mathematics in school (Leikin, 2018). Furthermore, Leikin points out that there is still a small number of studies in cognitive neuroscience that explore brain processing associated with more advanced mathematical concepts and that *these are rarely connected to theories in mathematics education*. The lack of applicability of cognitive neuroscience to the classroom is also stressed by Bowers (2016), who points out that there are no current examples of neuroscience motivating new and effective teaching methods, particularly due to the fact that characterizing children's cognitive capacities through behavioral measures is easier than using brain measures.

Nonetheless, Leikin argues that neurocognitive investigation can enrich mathematics education by helping to test conclusions drawn from behavioral research in mathematics education. For this purpose, neurocognitive research goals and questions can be informed by results from behavioral research, while behavioral studies can guide neuroscientific research through task design and interpretation of findings (Leikin, 2018).

Although this type of interdisciplinary research has the potential to influence the development of mathematics education, its significant applications in the classroom remain distant. Many programs that aim to bring neuroscience closer to teachers are based on neuromyths—for example, the belief that students learn better when information is presented in their preferred learning style (Privitera, 2021). Moreover, it is often teachers' own interest, combined with their lack of specialized training in neuroscience, that may partly contribute to the persistence of these neuromyths in education (Dekker et al., 2012). Thus, time and careful preparation are necessary before findings can be useful for teaching practice.

**The risk of narrowing the scope of research in mathematics education.**

A crucial question for our debate is what conditions must be met for research to inform teaching practice. From an organizational science perspective, direct knowledge exchange—without intermediaries—between research and practice is straightforward only when both domains share purposes, common understandings, and a low need for innovation (Carlile, 2004). In education, these conditions are not met: researchers and teachers pursue different objectives, reforms are constant, and teachers are neither dependent on nor typically resourced or motivated to apply evidence in their daily practice (Brown & Zhang, 2016). Therefore, translating research findings into the classroom demands considerable effort. Educators need to critically engage with research, focus on pressing questions, and reflect on their practice in light of empirical evidence. Moreover, practitioners have no opportunity to shape educational research, which limits its relevance for those responsible for implementing policy in teaching and learning contexts (Cain & Allan, 2017; Ball & Forzani, 2007).

Recognizing these challenges, Hiebert et al. (2002) propose building a professional knowledge base for teaching rooted in practitioner-generated insights—rich, context-sensitive findings capable of meaningfully impacting classroom instruction. They further advocate for a research-and-development system that integrates teachers' and researchers' expertise and unique skills, ensuring that practitioner knowledge is valued, shareable, and verifiable.

However, understanding teaching practice requires research questions to be framed and answered, which demands a theoretical perspective (Mason & Waywood, 1996; Hiebert & Stiegler, 2023). Theory provides the language to frame, describe, and discern phenomena. Thus, any research paradigm employs theory to formulate questions, specify methodology, define a coherent view of study objects, articulate consistent data collection and analysis approaches, and support the validation of findings (Mason & Waywood, 1996). Accordingly, fundamental research must not be underestimated. Advancing knowledge in mathematics education requires, for example, that researchers deepen theoretical concepts related to computational thinking, equity, and justice and develop new frameworks to study emerging classroom phenomena. Mathematics education—a discipline intimately linked to the vicissitudes of a changing society—demands that researchers continually question their prevailing methods and frameworks.

## Conclusion

Framing mathematics education research in terms of its classroom utility not only narrows the scope of what we understand as research in education but also shifts its value from the generation of knowledge to a mere instrument for achieving an external goal. However commendable that goal may be, this utilitarian approach contrasts the freedom of inquiry that characterizes research in other scientific disciplines.

While it is of utmost importance that research in mathematics education informs educational practice, it also serves broader purposes. For the field to thrive intellectually, it must preserve space for foundational, interdisciplinary, and disruptive inquiries—work that may not have immediate classroom applications but enriches the discipline in ways that may ultimately benefit educational practice.

## MATHEMATICS EDUCATION RESEARCH EMBEDED IN CLASSROOM PRACTICE (MAITREE INPRASITHA)

### The historical relationship between mathematics education and research

The history of mathematics education and mathematics education research reveals an important but often unfulfilled promise. Mathematics education as a field has always been concerned with a fundamental question: how can mathematics be taught and learned in ways that are meaningful to students? This question led to the emergence of mathematics education research as a systematic inquiry to inform classroom practices (Kilpatrick, 1992; Schoenfeld, 2000; Sierpinska & Kilpatrick, 1998). The very existence of our community is predicated on the belief that research should enhance mathematics teaching and learning. Looking at the historical development of mathematics education research, we can observe various paradigm shifts in approaches, methodologies, and focus areas.

The establishment of organizations like the International Commission on Mathematical Instruction (ICMI) in 1908 and its affiliated study groups, along with respected journals such as Educational Studies in Mathematics (first published in 1968) and the Journal for Research in Mathematics Education (1970), all share the stated aim of supporting mathematics teaching and learning. The International Group for the Psychology of Mathematics Education (PME), established in 1976, was founded precisely to bridge psychological research on mathematical thinking with pedagogical practices in classrooms (Kilpatrick, 1992). An examination of these journals' mission statements reveals consistent commitments to "linking theory and practice in mathematics education" (ESM) and "improving mathematics education through research" (JRME). The Mathematics Teacher Educator journal, established more recently in 2012, explicitly states its purpose as "contributing to the improvement of mathematics teacher education".

However, Lester (1994), in his comprehensive review in the Journal for Research in Mathematics Education, outlined how research paradigms have evolved over

decades—from behaviorist approaches through cognitive perspectives to sociocultural frameworks—each with their strengths and limitations. Yet despite these developments, a persistent critique has remained: much of this research fails to significantly improve mathematics teaching and learning in schools. Despite these laudable aims and the proliferation of research output over the past five decades, international assessments such as PISA (Programme for International Student Assessment) show that mathematics achievement worldwide has not significantly improved in the last 20 years. This suggests a fundamental disconnect between the growing body of mathematics education research and actual classroom practices. We must ask ourselves: Why has our research, despite its volume and methodological sophistication, not led to more substantial improvements in mathematics classrooms globally?

**The disconnect between research and classroom practice**

As a mathematics teacher in a Bangkok secondary school before becoming a mathematics teacher educator in 1987, I observed firsthand the gap between research and practice. This disconnects stems from several issues that our mathematics education community must address.

First, researchers often position themselves as external observers, adopting perspectives from social science that privilege theorizing over practical applications (Isoda, 2025, forthcoming). While theoretical frameworks are essential, they sometimes become ends in themselves rather than means to improve classroom practice. Mathematics education research has increasingly valued theoretical contributions and methodological innovations over tangible improvements in teaching and learning.

Second, the contexts in which research is conducted often differ significantly from the everyday realities of regular classrooms. Research settings frequently involve ideal conditions, special interventions, or highly motivated participants that do not reflect typical teaching situations. This creates an implementation gap where findings, while valid in research settings, prove difficult to translate into everyday practice. Guo et al. (2023) highlights the challenges of integrating simulation-based teaching with traditional methods and the tension between generative interactions and structured classroom culture. These contradictions illustrate the difficulties teachers face in applying research findings in practice, suggesting that effective co-design processes can help bridge this gap by fostering collaboration and perspective sharing between teachers and researchers.

Third, the language and communication of research findings are frequently inaccessible to practitioners. Research reports are written primarily for other researchers, using specialized terminology and focusing on issues that may not align with teachers' immediate concerns. This creates a two-community problem (Huberman, 1990) where researchers and practitioners operate in separate discourse communities with limited meaningful exchange.

Fourth, traditional research dissemination models assume a linear transfer of knowledge from researchers to practitioners, especially fundamental values whereby research merely *informs* practices. This model fundamentally misunderstands how teachers develop their practice, which typically occurs through collaborative reflection, experimentation, and adaptation within professional communities rather than through reading research reports.

As a result of these disconnects, we find ourselves in a situation where mathematics education research has grown substantially in volume and sophistication but has had limited impact on everyday classroom practices. This leads us to question not just the usefulness of our research but its very purpose and value.

## Japanese lesson study: research embedded in lived classrooms

One approach that bridges research and practice is Japanese Lesson Study, which gained international attention following Stigler and Hiebert's influential book *"The Teaching Gap"* (1999). Lesson Study represents a research approach deeply embedded in classroom practice, where teachers collaboratively plan, observe, and refine lessons (or classes) in a systematic cycle of inquiry. Many researchers employed Lesson Study as research methodology, such as Berkics (2024), who defines Lesson Study as a research methodology in educational settings by facilitating systematic, cooperative, and critical examination of teaching practices. It effectively addresses research questions while enhancing teaching practices, as demonstrated by its successful implementation in familiarizing students with model-based reasoning (Jansen et al., 2021). Arévalo et al. (2024) highlights Lesson Study as a valuable research methodology in educational settings, demonstrating its effectiveness in enhancing teachers' practical knowledge and understanding.

The global interest in Lesson Study has grown tremendously since 1999, with publications on the approach increasing from fewer than 10 English-language articles prior to 2000 to hundreds of publications annually in recent years. The establishment of International Journal for Lesson and Learning Studies (IJLLS) in 2012 is evident of the growing interest around the world. This growth also reflects recognition of Lesson Study's potential to connect research with improvements in classroom teaching.

What distinguishes Lesson Study is its situation within "lived classrooms"—real educational contexts where the complexities and challenges of everyday teaching are confronted directly. Unlike research approaches that abstract from classroom realities, Lesson Study embraces these complexities as essential to understanding and improving mathematics teaching and learning.

Key features that make Lesson Study valuable for connecting research and practice include:

1. Research is conducted *by* teachers, not just *on* teachers, repositioning practitioners as active knowledge producers rather than passive research subjects.

2. The research process is collaborative, involving multiple stakeholders including teachers, administrators, and university-based educators in a community of practice.
3. The focus remains firmly on student learning and thinking, grounding all inquiry in observable evidence from actual classroom interactions.
4. Knowledge is developed through iterative cycles that integrate planning, action, observation, and reflection, creating a continuous improvement model.
5. The research findings are immediately applicable, being developed in the same context where they will be used.

This approach represents a fundamentally different relationship between research and practice—one where the boundaries between researcher and practitioner blur, and where the ultimate value of research is judged by its contribution to improving teaching and learning. This is a paradigm shift from *research just informing practice* to *practice* (i.e., of Lesson Study) as *research.*

**The TLSOA model: transforming research in Thai classrooms**

Since 1999, I have been engaged in a longitudinal study transforming Japanese Lesson Study for use in Thai classrooms through what we have termed the Transformative Lesson Study incorporating Open Approach (TLSOA) model. This work has provided concrete evidence of how research can be genuinely useful for classroom practice when it is embedded within a school-based professional learning community.

The TLSOA model positions the school as a research site where multiple stakeholders—including teachers, administrators, university researchers, and teacher educators—collaborate to improve mathematics teaching and learning. The "Open Approach" serves as a thematic framework for implementing Lesson Study, emphasizing open-ended problems and multiple solution pathways. This approach has demonstrated significant benefits for both students and teachers. Students develop 21st-century thinking skills through engagement with open-ended mathematical problems and collaborative problem-solving. Teachers experience continuous professional development within a supportive community, developing deeper mathematical knowledge and pedagogical skills through collaborative planning, observation, and reflection (Inprasitha, 2022).

What makes this approach particularly valuable is that research occurs in authentic classroom settings, with researchers and practitioners working together on the "same front line" to understand and improve mathematics teaching and learning. The complexity of classroom practice is not reduced or simplified but becomes the central focus of collaborative inquiry.

The Open Class component of TLSOA creates a continuum for professional learning across multiple levels of the education system. Pre-service teachers, in-service teachers, school administrators, teacher educators, and researchers all participate in observing and analysing lessons, creating a shared space for developing practical knowledge about mathematics teaching and learning.

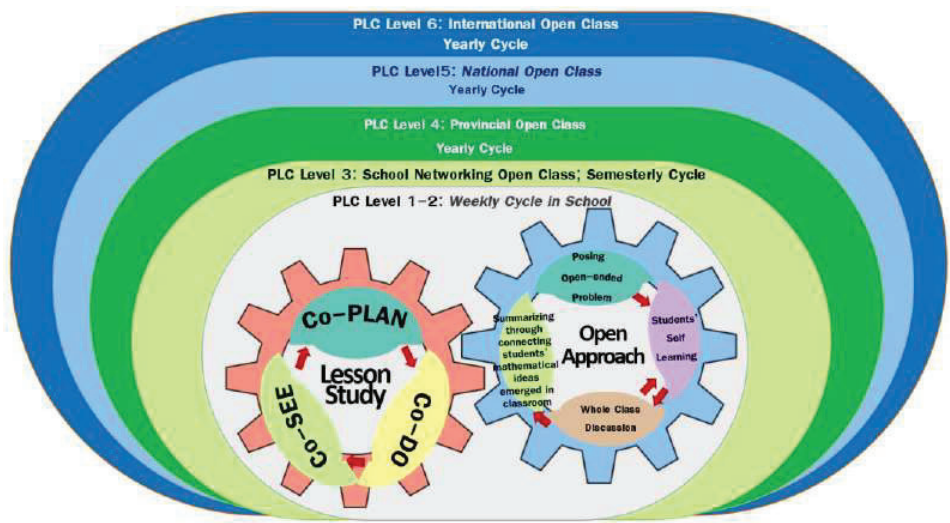

Figure 1: TLSOA Model: Open Class as a Continuum of Multi-Level Professional Learning Communities (adapted from Inprasitha, 2022)

Our longitudinal research has documented how this approach leads to sustainable improvements in classroom practice. Rather than research findings sitting unused in academic journals, knowledge generated through TLSOA is immediately implemented in teaching practices and continues to evolve through ongoing cycles of collaborative inquiry.

**Conclusion: towards research that serves the classroom**

Mathematics education research must be useful for the classroom—this is not merely a preference but an ethical imperative. If our research does not ultimately improve the teaching and learning experiences of students and teachers, we must question its purpose and value.

The history of mathematics education research shows that despite good intentions, much research has failed to significantly impact classroom practice. To address this disconnect, we need to reconceptualize the relationship between research and practice in several ways:

1. Research should be conducted *with* teachers rather than *on* teachers, recognizing practitioners as essential partners in knowledge generation.
2. The lived realities of classrooms must be the central context for research, embracing rather than simplifying the complexities of teaching and learning.
3. Knowledge development should occur within professional learning communities that span traditional boundaries between schools and universities.

4. The value of research should be judged primarily by its contribution to improving teaching and learning practices rather than by theoretical sophistication alone.

The Japanese Lesson Study approach and our TLSOA model in Thailand offer concrete examples of how research can be genuinely useful for classroom practice. These approaches demonstrate that when research is deeply embedded in the realities of classroom teaching, it can lead to sustainable improvements in mathematics education.

As a community of mathematics educators and researchers, we must commit to ensuring that our research serves its ultimate purpose: improving the mathematical experiences and outcomes of students in classrooms around the world. This requires not just new methodologies but a fundamental reorientation of our values and priorities as a field. Research that is not useful for the classroom is research that has forgotten its most fundamental purpose.

## MATHEMATICS EDUCATION RESEARCH MUST BE USEFUL FOR THE CLASSROOM (DEMETRA PITTA-PANTAZI)

During this panel discussion, I argue that Mathematics Education Research (MER) must be useful for the classroom—not merely as an academic pursuit but as a transformative force. To support this position, I briefly refer to the nature of research in mathematics education and the concept of today's classroom. Then, I present arguments about why research should be useful in the classroom.

### The nature of inquiry within mathematics education

MER has two core purposes: pure research, which is to understand mathematical thinking, teaching, and learning, and applied research, which is to use such understanding to improve mathematics instruction (Schoenfeld, 2003). It can be argued that prioritizing the applied research risks narrowing its intellectual scope, limiting its broader contributions to the field and society, and reducing it to immediate classroom utility. However, even the most 'pure' MER ultimately seeks to inform and enhance mathematical teaching and learning in classrooms. This objective underpins our scholarly endeavors across the diverse landscape of MER, establishing a unifying thread throughout seemingly disparate methodological approaches and theoretical frameworks.

Kilpatrick (2008) asserted that MER's theoretical and practical needs should be addressed scientifically rather than pragmatically to ensure that research remains useful for the classroom. Peter Taylor's question, "What's it all for?" (Kajander, 2020, p. 776), highlights the essential purpose of MER. If research is to fulfill its true purpose, it must provide tangible benefits for educators and students, ensuring classrooms become centers of effective and meaningful learning.

### The nature of the classroom

Classrooms are evolving rapidly due to social and technological changes. They are diverse and exist in various settings, such as homes, hospitals, museums, streets, and other informal settings. They may be in-person or online, synchronous or asynchronous. Classrooms involve multiple partners, including parents/guardians, informal schools, media creators, and textbook authors. The classroom model has changed and will change even more. Given this fluidity, MER must be useful for all these types of classrooms to address the various educational and global challenges.

### Why research must be useful for the classroom

Below, I highlight some arguments as to why mathematics education must be useful for the classroom. They are organized into two categories: academic and societal.

### Academic arguments

*Address critical educational needs*

MER must address the most pressing challenges faced by educators and students. Alignment with these needs ensures that research yields meaningful and applicable classroom outcomes. As Cai et al. (2017) stated, "If research is not solving teachers' problems, why would teachers want to use the research findings?" (p.345), highlighting the necessity for research to target areas requiring urgent improvements.

*Bridge the gap between theory and practice*

The traditional divide between research and practice has limited the impact of valuable findings on classroom instruction. MER gains credibility and advances rapidly when its theoretical constructs are regularly tested against the complex variables of actual classroom settings. MER must be useful for the classroom to bridge the gap between theoretical findings and teaching practices. Without a direct impact on practice, research remains an abstract exercise with little influence on student learning.

*Improve pedagogical strategies*

Research that directly informs practice helps teachers design and implement effective strategies. This symbiotic relationship between research and the classroom strengthens effective teaching strategies and, at the same time, improves research methodologies and theoretical frameworks, creating a continuous cycle of improvement.

*Enhance the dynamic exchange of knowledge and concerns between researchers and practitioners*

If MER is useful for the classroom, it will strengthen the collaboration between researchers and practitioners. Lasting improvements in education occur when practitioners actively define problems and develop context-specific solutions (Tyack & Cuban, 1995). Practitioners offer invaluable insights about their classrooms, while researchers contribute their expertise in literature, methodologies, and related topics. Strengthening the collaboration between practitioners and researchers enhances the understanding of the challenges that practitioners face (Cai et al., 2018). This collaboration grounds research in classroom realities and makes teaching more research-informed. Such synergies create a dynamic exchange of knowledge and concerns, allowing practitioners and researchers to refine their roles and enhance their impact.

*Facilitate complex decision-making in educational environments*

Educators work in complex, unpredictable classrooms requiring both theoretical knowledge and practical experience. Research provides the frameworks and evidence-based approaches to enhance decision-making and effectiveness. Referencing established findings helps practitioners justify pedagogical choices while preparing them for difficult or unusual situations by offering potential solutions.

*Encourage more teachers to participate in research*

Research designed for classroom utility has a broader impact on educational systems. When teachers see its benefits, they are more inclined to engage in research and encourage other teachers to do the same. This contributes to a shared body of knowledge that benefits the entire profession (Cai et al., 2018). The collective engagement creates a ripple effect, refining instructional practices across diverse educational settings and amplifying the research's impact on education.

*Address different students' needs at different locations*

Researchers can design appropriate instructional activities and adapt teaching to the needs of different students in different educational contexts. This adaptability ensures that research findings remain relevant and effective across different learning environments. Additionally, research useful for the classroom helps identify variables (e.g., social, personal, technological) that influence learning outcomes, leading to more accurate insights. In this way, research informs practice, provides valuable data on student learning, and leads to continuous advancements in students' learning.

*Identify potential issues before they become a widespread problem*

Research is useful not only because it sheds light on current problems but also helps detect emerging issues. By collecting and analyzing data across classrooms, researchers can pose and investigate new questions about teaching and learning. This forward-looking approach allows mathematics education to evolve and address emerging challenges. Identifying potential issues early enables research to serve a preventive function, helping educators develop proactive rather than reactive strategies.

*Create a cumulative knowledge base*

For research to have a lasting impact on classroom practice, it must contribute to an accessible, shared, and evolving knowledge base. As mathematics learning builds cumulatively, educational research should systematically build on previous findings rather than produce isolated studies. The NCTM (n.d) emphasizes that research-practice knowledge must be continuously updated with rigorous methods and reported in ways that can be shared, replicated, and adapted across educational settings. This ensures that valuable insights reach educators rather than being confined to academic journals. By promoting both theory and practice, research creates the conditions for ongoing progress in education. As findings accumulate, they enhance our understanding and promote current practice, theory building, and future investigations.

*Support teacher professional development*

MER supports teachers' professional development by enhancing teachers' knowledge for teaching, improving classroom practices, and fostering reflective teaching. It also promotes the integration of new technologies and innovative teaching methods.

## Societal arguments

The view that MER must be useful for the classroom extends beyond mere academic discussions, emphasizing the need to recognize its broader societal role. When research is disconnected from classroom realities, it may produce important but isolated theoretical frameworks that cannot address the challenges teachers and students face. Thus, I will point out that MER, useful for the classroom, can serve society in at least three key ways: (a) build students' resilience and address global issues, (b) promote equity, and (c) influence policies and public opinion.

*Build students' resilience and address global challenges*

In today's evolving educational landscape, particularly post-pandemic, research for the classroom is more urgent than ever. Mathematics education fosters complex reasoning skills addressing global challenges like climate change, public health, and economic and world crises. MER explores how classroom implementation equips learners with essential skills while promoting a more resilient educational system, helping students see mathematics not as an isolated academic subject but as a tool for understanding and impacting their communities.

Rapid technological and social changes demand flexible and forward-thinking education. Students need curricula and adaptable learning approaches emphasizing fundamental principles and transferable skills rather than rigid content that may become obsolete. Research-based teaching fosters more effective learning environments that meet current and future demands. It empowers individuals to engage meaningfully with data, recognize patterns, and make informed decisions.

*Promote equity*

MER can identify systemic barriers that impede students' access to quality mathematics education and propose solutions that promote equity (Cai et al., 2018). It explores how classroom implementation can promote a more equitable educational system appropriate for all students, regardless of their backgrounds. This approach helps tackle broader social issues such as race, gender, class, and identity. It transforms MER from a purely academic endeavor to a means for fostering social change.

*Influence policies and public opinion*

Research restricted to academic circles has limited classroom impact. Systematic research is useful for the classroom and provides valuable insights to help educational systems make informed decisions and changes. To be truly helpful, MER must influence the broader educational landscape, shaping policies related to curricula, educational practices, assessment, and resource distribution. The NCTM emphasizes the importance of communicating research to stakeholders to enact systematic, large-scale improvements. Effective MER that is useful for the classroom is more likely to reach a wider audience and influence public opinion.

## Conclusion

MER, which is useful for the classroom, provides multiple benefits. It offers theoretical foundations and evidence-based strategies that inform practice, while classroom experiences generate new insights and questions that drive further research. A commitment to make mathematics education useful for the classroom is essential for professional responsibility and ongoing educational and social improvement. By fostering MER that is useful for the classroom, we can enhance the quality of mathematics education, improve student learning outcomes, impact the whole educational system, and ultimately support a more resilient future for all.